\def\marker{\>\hbox{${\vcenter{\vbox{
    \hrule height 0.4pt\hbox{\vrule width 0.4pt height 6pt
    \kern6pt\vrule width 0.4pt}\hrule height 0.4pt}}}$}\>}
\def\gpic#1{#1
     \medskip\par\noindent{\centerline{\box\graph}} \medskip}
\documentclass[12pt]{article}
\usepackage{pstricks,pst-node}
\usepackage{amsfonts,amsthm,amsmath}
\newtheorem{theorem}{Theorem}
\def\flr#1{\left\lfloor #1 \right\rfloor}
\def\ceil#1{\left\lceil #1 \right\rceil}

\begin{document}
\author{
Daniel Cranston\thanks{University of Illinois, dcransto@uiuc.edu}
\and I. Hal Sudborough\thanks{University of Texas---Dallas, hal@utdallas.edu}
\and Douglas B. West\thanks{University of Illinois, west@math.uiuc.edu.
Work supported in part by the NSA under Award No.~MDA904-03-1-0037.}}

\title{Short Proofs for Cut-and-Paste Sorting of Permutations}
%\date{March 9, 2005}
\maketitle
\begin{abstract}
We consider the problem of determining the maximum number of moves required to
sort a permutation of $[n]$ using cut-and-paste operations, in which a segment
is cut out and then pasted into the remaining string, possibly reversed.
We give short proofs that every permutation of $[n]$ can be transformed to the
identity in at most $\flr{2n/3}$ such moves and that some permutations
require at least $\flr{n/2}$ moves.
\end{abstract}

\section{Introduction}
The problem of sorting a list of numbers is so fundamental that it has been
studied under many computational models.  Some of these models require sorting
``in place'', where elements are moved around within an array of fixed length
using various allowed operations.

A well-studied example is that of sorting using reversals of substrings.
This is motivated by applications in measuring the evolutionary distance
between genomes of different species \cite{GPS}.  The restricted case in
which the reversed substring must be an initial portion of the permutation is
the famous ``pancake problem'' \cite{GP,HS}.  One can also give each integer a
sign to denote the orientation of the gene \cite{KST,GPS,HP};
in the case of prefix reversal, this becomes the ``burnt pancake problem''
\cite{CB}.

Sorting by shifts has also been considered
\cite{EE,THRS}.  Shifting one block past another is a
reasonable operation when the permutation is stored using a linked list.  A
restriction of this operation is inserting the head element between two later
elements \cite{AW}.  Some work has also been done in which both reversals and
shifts are allowed \cite{WDM}.

In this note, we consider unsigned permutations and allow a more powerful
operation that incorporates both block reversal and block transposition.
A {\it (cut-and-paste) move} consists of cutting a substring out of the
permutation, possibly reversing it, and then pasting it back into the
permutation at any position.  The storage model that supports this conveniently
is a doubly-linked list.

We study permutations of $[n]$, where $[n]=\{1,\ldots,n\}$.  A natural
measure of how close a permutation is to the identity permutation is the number
of pairs $\{i,i+1\}$ that occur consecutively (in either order).  Call such
pairs {\it adjacencies}.  At most three adjacencies are created at each move
(involving the two ends of the moved substring and the element(s) that
formerly were next to the substring.  For $n\ge4$, there are permutations with
no adjacencies, but the identity permutation has $n+1$ adjacencies (counting
the front and back), so there are permutations of $[n]$ that require at least
$\ceil{(n+1)/3}$ moves to sort.

By a simple insertion sort, every permutation can be sorted in at most $n-1$
moves.  Indeed, since every list of $n$ distinct numbers has a monotone
sublist of length at least $\sqrt{n}$ \cite{ES}, we can insert the remaining
elements one at a time into a longest monotone sequence, reversing the full
list at the end if necessary, to sort in at most $n-\sqrt{n}+1$ moves.

Thus trivially we have $\ceil{(n+1)/3}\le f(n)\le n-\sqrt{n}+1$, where
$f(n)$ is the worst-case number of cut-and-paste moves needed to sort a
permutation of $[n]$.  In Section 2, we prove that $f(n)\ge\flr{n/2}$,
obtaining many permutations of $[n]$ that require at least $\flr{n/2}$ moves to
sort.  In Section 3, we prove that $f(n)\le\flr{2n/3}$, by presenting an
algorithm that sorts any permutation of $[n]$ using at most $\flr{2n/3}$
moves.

Our upper bound is weaker than that of Eriksson et. al~\cite{EE}; they proved
an upper bound of $\flr{(2n-2)/3}$ moves for the model that allows only the
weaker block transposition moves.  However, their proof is fairly lengthy and
involves considerable structural analysis.  Thus our contribution here is a
significantly shorter argument by a different method of proof.  We use a
``weight function'' argument: the maximum weight of a permutation of $[n]$ is
$2n/3$, and our algorithm reduces the weight by (on average) 1 unit per move.
Because our sorting operation is more powerful than those studied earlier, our
lower bound is a new result.

\section{The Lower Bound}

We write a permutation of $[n]$ as a list of numbers within brackets, without
commas; for example, $\pi = [\pi_1\, \cdots \, \pi_n]$.  For discussion of
the lower bound, we prepend $\pi_0=0$ and postpend $\pi_{n+1}=n+1$ to a
permutation $\pi$.  A move is performed using three (not necessarily distinct)
{\it cut point} indices $i$ such that $\pi_i$ and $\pi_{i+1}$ become separated
by the move.  These cut points partition $\pi$ into four disjoint substrings
and yield three possible moves: the string between the first two or last two
cut points is inserted at the third cut point, possibly reversed.
Given cut points $i,j,k$ with $0\le i\leq j\leq k\le n$, the results of the
three legal moves on $\pi$ are
$$
\begin{array}{cccc}
\left[ \pi_{0} \cdots \pi_{i} \right. & \pi_{j+1} \cdots \pi_{k} & \pi_{i+1}
\cdots \pi_{j} & \left. \pi_{k+1} \cdots \pi_{n+1}\right], \\
\left[ \pi_0 \cdots \pi_i \right. & \pi_{j+1} \cdots \pi_k & \pi_j \cdots
\pi_{i+1} & \pi_{k+1} \cdots \left. \pi_{n+1}\right], \\
\left[ \pi_0 \cdots \pi_i \right. & \pi_k \cdots \pi_{j+1} & \pi_{i+1}
\cdots \pi_j & \left. \pi_{k+1} \cdots \pi_{n+1}\right]. \\
\end{array}
$$
To reverse a string in place, let $j=k$ in the second form or $i=j$ in the
third form.

The trivial lower bound of $(n+1)/3$ was obtained by considering adjacencies.
To improve this bound, define a {\it parity adjacency} to be a pair of
consecutive values in $\pi$ having opposite parity.  With $0$ and $n+1$ fixed
at the ends, the identity permutation has $n+1$ parity adjacencies.  The key is
to show that each move increases the number of parity adjacencies by at most 2.

\begin{theorem}
$f(n)\geq \flr{n/2}$ for every positive integer $n$.
\end{theorem}

\begin{proof}
A permutation with all even values before all odd values has one parity
adjacency if $n$ is even, two if $n$ is odd.  Since the identity permutation
has $n+1$ parity adjacencies (counting the ends), it thus suffices to show that
each move increases the number of parity adjacencies by at most 2.

Consider a move that inserts $\pi_k\cdots \pi_l$ between $\pi_m$ and
$\pi_{m+1}$, reversed or not.  The three newly consecutive pairs of values
are $\{\pi_{k-1}\pi_{l+1},\pi_m\pi_k,\pi_l\pi_{m+1}\}$ if the block is not
reversed and $\{\pi_{k-1}\pi_{l+1},\pi_m\pi_l,\pi_k\pi_{m+1}\}$ if it is
reversed.  If the number of parity adjacencies increases by 3, then each new
pair is a new parity adjacency and no parity adjacencies were destroyed.

The six values in the new pairs and broken pairs form a cycle of pairs that
must alternate between opposite parity ($\not\equiv$) and equal parity
($\equiv$).  We show these requirements below, depending on whether the moving
block is reversed.
$$
\begin{array}{cc}
&\pi_{k-1}\not\equiv\pi_{l+1}\equiv\pi_l
\not\equiv\pi_{m+1}\equiv\pi_m\not\equiv\pi_k\equiv\pi_{k-1},\\
\text{or}&\pi_{k-1}\not\equiv\pi_{l+1}\equiv\pi_l
\not\equiv\pi_m\equiv\pi_{m+1}\not\equiv\pi_k\equiv\pi_{k-1}.\\
\end{array}
$$

In each case, the six elements appear on a cycle whose steps alternate
between preserving parity and switching parity.  However, parity cannot
change an odd number of times along a cycle.  This contradiction prevents the
number of parity adjacencies from increasing by 3.
\end{proof}

The proof in fact shows that every input permutation with one parity
adjacency (when $n$ is even) or two parity adjacencies (when $n$ is odd)
requires $\flr{n/2}$ moves to sort.

\section{Upper Bound}

For the upper bound, we no longer append $0$ and $n+1$ at the ends of the
permutation.  Instead, we adopt the convention that the values $n$ and $1$ are
consecutive.  That is, adding 1 to the value $n$ produces the value $1$; all
computations with values are modulo $n$.

A {\it block} in a permutation is a maximal substring of (at least two)
consecutive values in consecutive positions.  A block is {\it increasing} if
its second value is one more than its first.  An increasing block can reach
$n$ and continue with 1.  A permutation that consists of a single increasing
block can be sorted in one move, shifting the initial part ending at $n$ to
the end.

An element in no block is a {\it singleton}.  When singletons $i$ and $i+1$
are made consecutive, they automatically form an adjacency.  When blocks ending
at $i$ and $i+1$ are made consecutive, they produce a larger block only if
oriented properly.  Hence somehow moves that reduce the number of blocks are
more valuable than moves that combine singletons.

We capture this phenomenon by giving higher weight to blocks than to singletons
in measuring the ``non-sortedness'' of a permutation.  The proof of the upper
bound is constructive, providing an algorithm to sort permutations of $[n]$
with at most $\flr{2n/3}$ moves.  We note that the algorithm for the pancake
problem in \cite{GP} also treats blocks and singletons differently, but not via
a weight function.

\begin{theorem}
$f(n)\leq \flr{2n/3}$, for every positive integer $n$.
\end{theorem}

\begin{proof}
We first give an algorithm that sorts the permutation using at most $\flr{2n/3} +1$ moves, then give a small modification that sorts the permutation in at most $\flr{2n/3}$ moves.

Let the {\it weight} of a permutation be
$$(number~of~blocks)+(2/3)(number~of~singletons).$$
Every permutation consisting of a single increasing block, including the
identity permutation, has weight 1.  The maximum weight of a permutation of
$[n]$ is $2n/3$, achieved by permutations having no blocks.  The {\it gain} of
a move is the amount by which it reduces the weight.

Using one move, we establish an increasing block at the beginning of
the permutation.  We show that when such a block exists, we can gain at
least 1 in one move or gain at least 2 in two moves, while maintaining that
condition.  We thus reach a permutation with weight 1 (a single increasing
block) in at most $1+\flr{2n/3}-1$ steps.  One more move may be needed
to sort the single increasing block.

In each move, the string we cut out is a union of full blocks and singletons,
never part of a block.  Also, we never paste this string into the interior of
a block.  That is, we never break adjacencies, and the numbers of singletons
and blocks can change only by creating adjacencies.  Creating a block from two
singletons gains $1/3$.  Absorbing a singleton into an existing block gains
$2/3$ (this is an {\it absorbing move}).  Creating one block by combining two
blocks gains $1$.  In all moves we ignore the possible gain resulting from
closing the gap left by the extracted string, counting only the gain from
pasting it elsewhere.

If our permutation has $i$ and $i+1$ in distinct blocks, then cutting out one
of these blocks and pasting it next to the other gains at least 1.  We call
this a {\it block move}.  When a block move is available, we perform it,
leaving (or augmenting) the increasing block at the beginning of the
permutation.

A {\it bonus move} is one that gains two adjacencies by the insertion.  That
is, the cut string (without splitting a block) has $i$ and $j$ at its ends,
and it is inserted between two consecutive elements whose values are next to
$i$ and next to $j$.  If one of these four elements was in a block, then the
move gains 1.  If at least two of them were in blocks, then the gain is at
least $4/3$ (this is an {\it extra bonus move}).  When we make two moves to
gain 2, they will be an absorbing move and an extra bonus move.

Let $l$ be the last value in the initial increasing block, and let $l'=l+1$.
Let $p$ be the value in the position following $l$, and let $p'$ be a value
immediately above or below $p$ that is not located next to $p$ (there are two
choices for $p'$ if $p$ is singleton, one if $p$ is in a block).  We consider
several cases based on the condition of these elements.

If $l'$ is in a block, then a block move is available, so we may assume that
$l'$ is a singleton.  If $p'$ is a singleton, or if the string $S$ with $l'$
and $p'$ as its end-elements contains the block that $p'$ ends, then we cut $S$
and insert it between $l$ and $p$ (if $p'$ is at the left end of $S$, we
reverse $S$ before inserting); see Figure 1.  Since $l$ is in a block, this
is a bonus move that gains at least 1 (extra bonus if $p'$ is in a block).

\gpic{
\expandafter\ifx\csname graph\endcsname\relax \csname newbox\endcsname\graph\fi
\expandafter\ifx\csname graphtemp\endcsname\relax \csname newdimen\endcsname\graphtemp\fi
\setbox\graph=\vtop{\vskip 0pt\hbox{%
    \graphtemp=.5ex\advance\graphtemp by 0.467in
    \rlap{\kern 0.667in\lower\graphtemp\hbox to 0pt{\hss $l$\hss}}%
    \graphtemp=.5ex\advance\graphtemp by 0.467in
    \rlap{\kern 1.111in\lower\graphtemp\hbox to 0pt{\hss $p$\hss}}%
    \graphtemp=.5ex\advance\graphtemp by 0.467in
    \rlap{\kern 1.778in\lower\graphtemp\hbox to 0pt{\hss $l'$\hss}}%
    \graphtemp=.5ex\advance\graphtemp by 0.467in
    \rlap{\kern 2.889in\lower\graphtemp\hbox to 0pt{\hss $p'$\hss}}%
    \special{pn 8}%
    \special{pa 1667 444}%
    \special{pa 1667 222}%
    \special{pa 3000 222}%
    \special{pa 3000 444}%
    \special{fp}%
    \special{pa 778 333}%
    \special{pa 0 333}%
    \special{pa 0 556}%
    \special{pa 778 556}%
    \special{pa 778 333}%
    \special{fp}%
    \special{pa 2333 222}%
    \special{pa 2333 89}%
    \special{fp}%
    \special{sh 1.000}%
    \special{pa 2311 178}%
    \special{pa 2333 89}%
    \special{pa 2356 178}%
    \special{pa 2311 178}%
    \special{fp}%
    \special{pa 2333 89}%
    \special{pa 2333 0}%
    \special{fp}%
    \special{pa 2333 0}%
    \special{pa 1467 0}%
    \special{fp}%
    \special{sh 1.000}%
    \special{pa 1556 22}%
    \special{pa 1467 0}%
    \special{pa 1556 -21}%
    \special{pa 1556 22}%
    \special{fp}%
    \special{pa 1467 0}%
    \special{pa 889 0}%
    \special{fp}%
    \special{pa 889 0}%
    \special{pa 889 233}%
    \special{fp}%
    \special{sh 1.000}%
    \special{pa 911 144}%
    \special{pa 889 233}%
    \special{pa 867 144}%
    \special{pa 911 144}%
    \special{fp}%
    \special{pa 889 233}%
    \special{pa 889 389}%
    \special{fp}%
    \graphtemp=.5ex\advance\graphtemp by 0.889in
    \rlap{\kern 1.422in\lower\graphtemp\hbox to 0pt{\hss Figure 1: A bonus move\hss}}%
    \hbox{\vrule depth0.978in width0pt height 0pt}%
    \kern 3.000in
  }%
}%
}

Hence we may assume that $p'$ is in a block $B$ that is not contained in $S$.
Let $q$ be the end of $B$ other than $p'$.  Let $q'$ be the value next to $q$
(up or down) that is not in $B$.  If $q'$ or $p$ is in a block, then a block
move is available to combine $B$ with the block containing $q'$ or $p$,
so we may assume that $q'$ and $p$ are singletons.  See Figure 2.

\gpic{
\expandafter\ifx\csname graph\endcsname\relax \csname newbox\endcsname\graph\fi
\expandafter\ifx\csname graphtemp\endcsname\relax \csname newdimen\endcsname\graphtemp\fi
\setbox\graph=\vtop{\vskip 0pt\hbox{%
    \graphtemp=.5ex\advance\graphtemp by 0.606in
    \rlap{\kern 2.788in\lower\graphtemp\hbox to 0pt{\hss $B$\hss}}%
    \graphtemp=.5ex\advance\graphtemp by 0.606in
    \rlap{\kern 0.727in\lower\graphtemp\hbox to 0pt{\hss $l$\hss}}%
    \graphtemp=.5ex\advance\graphtemp by 0.606in
    \rlap{\kern 1.212in\lower\graphtemp\hbox to 0pt{\hss $p$\hss}}%
    \graphtemp=.5ex\advance\graphtemp by 0.606in
    \rlap{\kern 1.939in\lower\graphtemp\hbox to 0pt{\hss $q'$\hss}}%
    \graphtemp=.5ex\advance\graphtemp by 0.606in
    \rlap{\kern 2.424in\lower\graphtemp\hbox to 0pt{\hss $q$\hss}}%
    \graphtemp=.5ex\advance\graphtemp by 0.606in
    \rlap{\kern 3.152in\lower\graphtemp\hbox to 0pt{\hss $p'$\hss}}%
    \graphtemp=.5ex\advance\graphtemp by 0.606in
    \rlap{\kern 3.879in\lower\graphtemp\hbox to 0pt{\hss $l'$\hss}}%
    \special{pn 8}%
    \special{pa 848 461}%
    \special{pa 0 461}%
    \special{pa 0 703}%
    \special{pa 848 703}%
    \special{pa 848 461}%
    \special{fp}%
    \special{pa 3273 461}%
    \special{pa 2303 461}%
    \special{pa 2303 703}%
    \special{pa 3273 703}%
    \special{pa 3273 461}%
    \special{fp}%
    \special{pa 1818 582}%
    \special{pa 1818 339}%
    \special{pa 4000 339}%
    \special{pa 4000 582}%
    \special{fp}%
    \special{pa 2788 703}%
    \special{pa 2788 921}%
    \special{fp}%
    \special{sh 1.000}%
    \special{pa 2812 824}%
    \special{pa 2788 921}%
    \special{pa 2764 824}%
    \special{pa 2812 824}%
    \special{fp}%
    \special{pa 2788 921}%
    \special{pa 2788 1067}%
    \special{fp}%
    \special{pa 2788 1067}%
    \special{pa 1755 1067}%
    \special{fp}%
    \special{sh 1.000}%
    \special{pa 1852 1091}%
    \special{pa 1755 1067}%
    \special{pa 1852 1042}%
    \special{pa 1852 1091}%
    \special{fp}%
    \special{pa 1755 1067}%
    \special{pa 1067 1067}%
    \special{fp}%
    \special{pa 1067 1067}%
    \special{pa 1067 812}%
    \special{fp}%
    \special{sh 1.000}%
    \special{pa 1042 909}%
    \special{pa 1067 812}%
    \special{pa 1091 909}%
    \special{pa 1042 909}%
    \special{fp}%
    \special{pa 1067 812}%
    \special{pa 1067 642}%
    \special{fp}%
    \special{pa 2788 339}%
    \special{pa 2788 194}%
    \special{fp}%
    \special{sh 1.000}%
    \special{pa 2764 291}%
    \special{pa 2788 194}%
    \special{pa 2812 291}%
    \special{pa 2764 291}%
    \special{fp}%
    \special{pa 2788 194}%
    \special{pa 2788 97}%
    \special{fp}%
    \special{pa 2788 97}%
    \special{pa 1697 97}%
    \special{fp}%
    \special{sh 1.000}%
    \special{pa 1794 121}%
    \special{pa 1697 97}%
    \special{pa 1794 73}%
    \special{pa 1794 121}%
    \special{fp}%
    \special{pa 1697 97}%
    \special{pa 970 97}%
    \special{fp}%
    \special{pa 970 97}%
    \special{pa 970 352}%
    \special{fp}%
    \special{sh 1.000}%
    \special{pa 994 255}%
    \special{pa 970 352}%
    \special{pa 945 255}%
    \special{pa 994 255}%
    \special{fp}%
    \special{pa 970 352}%
    \special{pa 970 521}%
    \special{fp}%
    \graphtemp=.5ex\advance\graphtemp by 1.164in
    \rlap{\kern 1.697in\lower\graphtemp\hbox to 0pt{\hss Move 1\hss}}%
    \graphtemp=.5ex\advance\graphtemp by 0.000in
    \rlap{\kern 1.697in\lower\graphtemp\hbox to 0pt{\hss Move 2\hss}}%
    \graphtemp=.5ex\advance\graphtemp by 1.552in
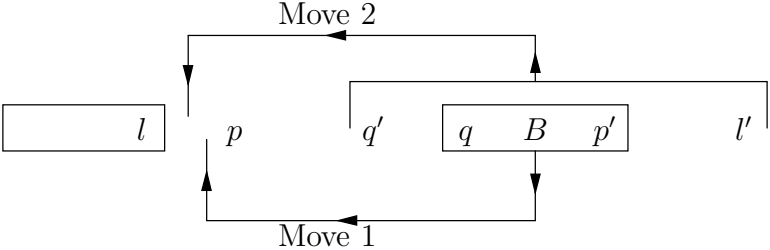
    \rlap{\kern 2.036in\lower\graphtemp\hbox to 0pt{\hss Figure 2: An absorbing move followed by an extra bonus move\hss}}%
    \hbox{\vrule depth1.648in width0pt height 0pt}%
    \kern 4.000in
  }%
}%
}

Now we make two moves.  We first insert $B$ between $l$ and $p$, with $p'$
next to $p$ (if $p'$ is at the left end of $B$, then we reverse $B$ before
inserting); this is an absorbing move, and it puts $q$ next to $l$.  We now
insert the string from $l'$ to $q'$ between $l$ and $q$ (if $q'$ is at the
left end of this string, we reverse this string before inserting); since both
$l$ and $q$ are now in blocks, this is an extra bonus move.  As noted earlier,
an absorbing move and an extra bonus move together gain (at least) 2.

A slight variation of the algorithm sorts the permutation in at most
$\flr{2n/3}$ moves.  The idea is to keep element $1$ in the first position and
thus avoid a final move to sort an increasing block.  In general, when the first
element is $a$, we bring a segment from $1$ to $a'$ to the front.  Now the
permutation begins $1b$, and we put a segment from $2$ to $b'$ between them
without breaking a block.  Six singletons become three blocks, or five become
two blocks (if $b=a'$), so the resulting permutation has weight at most
$2n/3-1$.  Since we have spent two moves and reduce to weight 1 at the end, we
use at most $\flr{2n/3}$ moves.

This argument may fail when $a\in \{1,2\}$ or $b=n$.  When $a=1$, we sort the
rest of the list inductively.  When $a=2$, we bring $1$ to the front and sort
the rest inductively in at most $2(n-2)/3$ moves.  When $b=n$, we first move
$1n$ to the front and then move the segment from $n$ till just before $2$ to
the back (reversed).  The permutation of $\{3,\ldots,n-1\}$ can now be sorted
inductively in $\flr{2(n-3)/3}$ moves, which with the two initial moves
satisfies the bound.
\end{proof}

It is straightforward to implement this algorithm with a running time of
$O(n^2)$.  At each move, we linearly search for each of the numbers indicated
in Figure 1 or Figure 2.  Jeff Erickson has noted that the algorithm can be
implemented with a running time of $O(n\log n)$, but the details are much
more complicated.

%\bibliographystyle{plain}
%\bibliography{pancakes}

\end{document}